\documentclass[10pt,a4paper,reqno]{amsart}
\usepackage{amsmath}
\usepackage{amssymb,latexsym}
\usepackage[dvips]{graphicx}
\usepackage{color}
\usepackage{cite}
\usepackage{amsthm}

\newcommand{\R}{\mathbb R}

\newcommand{\N}{\mathbb N}

\newtheorem{theorem}{Theorem} [section]
\newtheorem{lemma}{Lemma} [section]

\let\ssection=\section\renewcommand{\section}{\setcounter{equation}{0}\ssection}
\begin{document}
\title{ON THE GENERALIZED NONLINEAR CAMASSA-HOLM
EQUATION.}
\author{Mohamad Darwich, Samer Israwi and Raafat Talhouk}
\address{Department of Mathematics, Faculty of Sciences 1 and Laboratory of Mathematics, Doctoral School of Sciences and Technology, Lebanese University Hadat, Lebanon.}

\maketitle
\begin{abstract}
In this paper, a generalized nonlinear Camassa-Holm equation with time- and space-dependent coefficients is considered. We show that the control of the higher order dispersive  term is possible by using an adequate weight function to define the energy. The existence and uniqueness of solutions are obtained  via a Picard iterative method.
\end{abstract}


\section{Introduction}
\subsection{Presentation of the problem}
In this paper, we study  the Cauchy problem for the general nonlinear  higher order Camassa-Holm-type equation:
\begin{equation}\label{hch}
\left\{
	\begin{array}{l}
      (1-m\partial_x^2)u_t + a_1(t,x,u)u_x +a_2(t,x,u,u_x)u_{xx}
      \\\\\qquad+a_3(t,x,u) u_{xxx} +a_4(t,x) u_{xxxx}+a_5(t,x) u_{xxxxx}=f \quad\hbox{ for} \quad(t,x) \in(0,T]\times\R
	\\\\
	u_{\vert_{t=0}}=u^0,
	\end{array}\right.
\end{equation}
where $u=u(t,x)$, from $[0,T]\times\R$ into $\R$, is the unknown function of the problem, 
$m>0$ and $a_i$, $1 \leq i \leq 5$, are real-valued smooth given
functions  where their exact regularities will be precised later. This equation covers several important unidirectional models for the water waves problems at different regimes which take into account the  variations of the bottom. We have in view in particular the  example of the Camassa-Holm  equation (see\cite{m2an}), which is more nonlinear then the  KdV equation (see for instance \cite{CRASS},\cite{kdv from vincent}, \cite{craigkdv}, \cite{Linaresponce}, \cite{TianGao}). However, the most prominent example that we have in mind is the Kawahara-type  approximation ( see \cite{Erika}), in which case the coefficient $a_5$ does not vanish. The presence of the fifth order derivative term is very important, so that the  equation describes both  nonlinear and dispersive effects as does the Camassa-Holm equation in the case of special tension surface values (see \cite{lannes1}).\\
 Looking for solutions of  (\ref{hch}) plays an important and significant role in the study of unidirectional limits for water wave problems with variable depth and topographies. To our knowledge the problem (\ref{hch}) has not  been analyzed previously. In the present paper, we prove the local  well-posedness of the initial value problem (\ref{hch})  by a 
standard Picard iterative scheme and the use of adequate energy estimates under a condition of nondegeneracy of the higher dispersive  coefficient $a_5$.  

\subsection{Notations and Main result}

 In the following, $C_0$ denotes any nonnegative constant whose exact expression is of no importance. The notation $a\lesssim b$ means that 
 $a\leq C_0\ b$.\\
 We denote by $C(\lambda_1, \lambda_2,\dots)$ a nonnegative constant depending on the parameters
 $\lambda_1$, $\lambda_2$,\dots and whose dependence on the $\lambda_j$ is always assumed to be nondecreasing.\\
 For any $s \in \R$, we denote $[s]$ the integer part of $s$.\\
 Let $p$ be any constant
 with $1\leq p< \infty$ and denote $L^p=L^p(\R)$ the space of all Lebesgue-measurable functions
 $f$ with the standard norm \[\vert f \vert_{L^p}=\big(\int_{\R}\vert f(x)\vert^p dx\big)^{1/p}<\infty.\] The real inner product of any two functions $f_1$
 and $f_2$ in the Hilbert space $L^2(\R)$ is denoted by
\[
 (f_1,f_2)=\int_{\R}f_1(x)f_2(x) dx.
 \]
 The space $L^\infty=L^\infty(\R)$ consists of all essentially bounded and Lebesgue-measurable functions
 $f$ with the norm
\[
 \vert f\vert_{L^\infty}= \sup \vert f(x)\vert<\infty.
\]
 We denote by $W^{1,\infty}(\R)=\{f, \mbox{ s.t. }f,\partial_x f\in L^{\infty}(\R)\}$ endowed with its canonical norm.\\
 For any real constant $s\geq0$, $H^s=H^s(\R)$ denotes the Sobolev space of all tempered
 distributions $f$ with the norm $\vert f\vert_{H^s}=\vert \Lambda^s f\vert_{L^2} < \infty$, where $\Lambda$
 is the pseudo-differential operator $\Lambda=(1-\partial_x^2)^{1/2}$.\\
 For any two functions $u=u(t,x)$ and $v(t,x)$
 defined on $ [0,T)\times \R$ with $T>0$, we denote the inner product, the $L^p$-norm and especially
 the $L^2$-norm, as well as the Sobolev norm,
 with respect to the spatial variable $x$, by $(u,v)=(u(t,\cdot),v(t,\cdot))$, $\vert u \vert_{L^p}=\vert u(t,\cdot)\vert_{L^p}$,
 $\vert u \vert_{L^2}=\vert u(t,\cdot)\vert_{L^2}$ , and $ \vert u \vert_{H^s}=\vert u(t,\cdot)\vert_{H^s}$, respectively.\\
 We denote $L^\infty([0,T);H^s(\R))$ the space of functions such that $u(t,\cdot)$ is controlled in $H^s$, uniformly for $t\in[0,T)$: $\big\Vert u\big\Vert_{L^\infty([0,T);H^s(\R))} \ = \ \sup_{t\in[0,T)}\vert u(t,\cdot)\vert_{H^s} \ < \ \infty.$\\
 Finally, $C^k(\R^i)$, $i\ge 1$ denote the space of
 $k$-times continuously differentiable functions.\\
 For any closed operator $T$ defined on a Banach space $X$ of functions, the commutator $[T,f]$ is defined
  by $[T,f]g=T(fg)-fT(g)$ with $f$, $g$ and $fg$ belonging to the domain of $T$. The same notation is used for $f$ as an operator mapping the domain of $T$ into itself.\\
 Moreover, we define the following operators:
  $\displaystyle{\Lambda_m^s=(1-m\partial_x^2)^{\frac{s}{2}}}$ and its inverse $\displaystyle{\Lambda_m^{-s}}$ such that $\displaystyle{\displaystyle{\hat{\Lambda}_m^{-s}(u)}= (1+m\xi^2)^{-\frac{s}{2}}\hat{u}}.$
  
  Finally, we will study the local well-posedness of the  initial value problem (\ref{hch}) in $H^{s}(\R)$  endowed with canonical norm.\\
Let us now state our main result:
\begin{theorem}\label{th1}
Let $s>\frac{5}{2}$  and $f\in C([0,T];H^{s}(\R))$.  We suppose that: 
	\begin{itemize}
	
	\item $a_1,a_3$  in $C([0,T],C^{[s]+1}(\R^2))$, $a_2$ in $C([0,T],C^{[s]+1}(\R^3))$ and $\partial ^k_{x} a_j$ are bounded with respect to $x$ for all $0 \leq k \leq [s] +1$ and $ 1 \leq j \leq 3$.
         \item $a_4\in C([0,T];H^{s+1}(\R))$, $\partial_t a_4 \in L^{\infty}(0,T,L^{\infty}(\mathbb{R}))$
         \item $a_5\in C([0,T];H^{s+2}(\R))\cap C([0,T],L^{\infty}(\R))$ \hbox{ with} $ \partial_ta_5\in L^{\infty}(0,T;L^{\infty}(\R))$,
         \item $F(t,x):=\int_0^x\frac{a_4}{a_5}dy \in C([0,T];L^{\infty}(\R)) ~\text{and} ~\partial_tF \in L^{\infty}(0,T;L^{\infty}(\R))$,

         \end{itemize}
         Assume moreover that 
         there is a positive constant $ c_1>0 $ such that $c_1\leq\vert a_5(t,x)\vert \;\forall \,(t,x) \in [0,T]\times\R.$ 
  	Then for all $u^0 \in H^{s}(\R)$, there exist a time $T^{\star}>0$ and a unique solution 
	$u$ to (\ref{hch})
	in $C([0,T^{\star}];H^{s})$.
\end{theorem}
  
\section{Proof of the Main results}
Before we start the proof, we  give the following useful lemma:
\begin{lemma}
 Let $m>0$, $s\in \R^+$  then the linear operator $\Lambda^2_m$: $H^{s+2}(\R)\to H^s(\R)$  is well defined, continuous, one-to-one and onto.
  If we suppose that $u=\Lambda^{-2}_mf $ for $f\in H^s(\R)$ then\\
  \begin{eqnarray}\label{est1}
\vert u\vert_{H^{s+2}}\ & \leq& \  \frac{1}{m}\vert f\vert_{H^s}\quad \hbox{if}\quad 0<m\leq1 \\
\vert u\vert_{H^{s+2}}\ & \leq& \  \vert f\vert_{H^s} \quad \hbox{if}\quad m\ge 1.
\end{eqnarray}
Moreover,$$\Lambda^s\Lambda^{-2}_m=\Lambda^{s-2}\Lambda^{0}_m=\Lambda^{0}_m\Lambda^{s-2},$$ 
where 
$\Lambda^0_m$: $H^{s}(\R)\to H^s(\R)$ is linear continuous  one-to-one and onto operator defined by 

$$\displaystyle{\widehat{\Lambda^0_m} u(\xi)=(1+\xi^2)(1+m\xi^2)^{-1}\hat{u}(\xi),}$$

 with 
\begin{eqnarray}\label{est2}
\vert \Lambda^0_m \vert_{H^{s}\to H^s}\ & \leq& \  \max{(\frac{1}{m},1)},\\
\vert (\Lambda^0_m)^{-1} \vert_{H^{s}\to H^s}\ & \leq& \  \max{(m,1)}.
\end{eqnarray}
  \end{lemma}
  \proof $\|\Lambda^{-2}_mf\|_{H^{s+2}} = \|(1+\xi^2)^{\frac{s}{2}+1}(1+m\xi^2)^{-1}\hat{f}\|_{L^2}$.\\
  If $ m \geq 1$, then $1+m\xi^2 \geq 1+ \xi^2$ and $\frac{ 1+ \xi^2}{ 1+ m\xi^2} \leq 1$, then $\|(1+\xi^2)^{\frac{s}{2}+1}(1+m\xi^2)^{-1}\hat{f}\|_{L^2}= \|(1+\xi^2)^{\frac{s}{2}}(1+\xi^2)(1+m\xi^2)^{-1}\hat{f}\|_{L^2} \leq \|(1+\xi^2)^{\frac{s}{2}}\hat{f}\|_{L^2}$.\\
  If $0< m < 1 $, we have $\frac{ 1+ \xi^2}{ 1+ m\xi^2}= 1 + (1-m)\frac{\xi^2}{ 1+ m\xi^2} \leq 1 + \frac{(1-m)}{m} = \frac{1}{m}$, then\\ 
  $\|\Lambda^{-2}_mf\|_{H^{s+2}} \leq \frac{1}{m} \|f\|_{H^{s}}$.\\
  Now $\|\Lambda^{0}_mf\|_{H^s} = \|\Lambda^2\Lambda^{-2}_mf\|_{H^s}=\|\Lambda^{-2}_mf\|_{H^{s+2}} \leq \max(1,\frac{1}{m})\|f\|_{H^s}$.\\
$\|(\Lambda^{0}_m)^{-1}f\|_{H^s} = \|(1+m\xi^2)(1+\xi^2)^{-1}(1+\xi^2)^{{\frac{s}{2}}}\hat{f}\|_{L^2}$.\\
  If $ m \geq 1$, then $(1+m\xi^2)(1+\xi^2)^{-1} = 1 + (m-1)\frac{\xi^2}{1+\xi^2} \leq m$, then $\|(\Lambda^{0}_m)^{-1}f\|_{H^s} \leq m \|f\|_{H^s}.$\\
  If $0< m < 1$, $(1+m\xi^2)(1+\xi^2)^{-1} \leq 1$, then $\|(\Lambda^{0}_m)^{-1}f\|_{H^s} \leq \|f\|_{H^s}$.\\Finally $\|(\Lambda^{0}_m)^{-1}f\|_{H^s} \leq \max(1,m) \|f\|_{H^s}$.\\ 
  \\
   We will start the proof of Theorem \ref{th1} by studying a linearized problem associated to (\ref{hch}):
\subsection{Linear analysis:}
For any smooth enough $v$, we define the ``linearized'' operator:
\begin{eqnarray*}
	{\mathcal L}(v,\partial)&=&\Lambda^2_m\partial_t +a_1(t,x,v)\partial_x+a_2(t,x,v,v_x)\partial_x^2 + 
        a_3(t,x,v) \partial_x^3+a_4(t,x) \partial_x^4+a_5 (t,x)\partial_x^5
       \end{eqnarray*}
and  the following  initial value problem: 
\begin{equation}\label{lhch}
	\left\lbrace
	\begin{array}{l}
	{\mathcal L}(v,\partial)u=f,\\
	u_{\vert_{t=0}}=u^0.
	\end{array}\right.
\end{equation}

Equation (\ref{lhch}) is a linear equation  which can be solved by a standard method (see \cite{Taylor}) in any time interval in which its coefficients are defined and regular enough. 
We first establish some precise energy-type estimates of the solution. 
We define the ``energy'' norm, 
$$
	E^s(u)^2=\vert w \Lambda^s u\vert_{L^2}^2,
$$
where $w$ is a weight function that will be chosen later. 
For the moment, we just require that 
 there exists  two positive numbers $w_1,w_2$
such that for all $(t,x)$ in $(0,T]\times\R$,
$$
w_1\leq w(t,x)\leq w_2,
$$
so that  $E^s(u)$ is uniformly equivalent to the standard $H^s$-norm.
Differentiating $\frac{1}{2}e^{-\lambda t}E^s(u)^2$
with respect to time, one gets using  (\ref{lhch})
\begin{eqnarray*}\label{method}
	 \lefteqn{\frac{1}{2}e^{\lambda t}\partial_t(e^{-\lambda t} E^s(u)^2)
	=-\frac{\lambda}{2} E^s(u)^2 -
         \big(\Lambda_m^0\Lambda^{s-2} (a_1u_x) ,w^2\Lambda^s u\big)}\\\nonumber
	&&-\big(\Lambda_m^0\Lambda^{s-2}( a_2u_{xx}) ,w^2\Lambda^su\big)-\big(\Lambda_m^0\Lambda^{s-2}(a_3 u_{xxx}),w^2\Lambda^s u\big)
	-\big(\Lambda_m^0\Lambda^{s-2}( a_4u_{xxxx}),w^2\Lambda^s u\big)\\\nonumber
	& &- \big(\Lambda_m^0\Lambda^{s-2}(a_5u_{xxxxx}),w^2\Lambda^s  u\big)+\big(\Lambda_m^0\Lambda^{s-2}f ,w^2\Lambda^su\big)
        +\big(2ww_t\Lambda^su ,\Lambda^su\big).\\\nonumber
\end{eqnarray*}
We now turn to estimating the different terms of the r.h.s  of the previous identity. 
\\
$\bullet$ Estimate of $\big(\Lambda^{s-2} (a_1u_x ),\Lambda^0_m w^2\Lambda^su\big)$. By the Cauchy-Schwarz inequality we have  \begin{align}
\vert \big(\Lambda^{s-2} (a_1u_x ),\Lambda^0_m w^2\Lambda^{s}u\big)\vert&\leq \frac{1}{m}\vert a_1(t,x,v)\vert_{H^{s-2}}\vert u_x \vert_{H^{s-1}} \vert w^2\Lambda^s u \vert_{L^2}\nonumber \\ &\leq 
C(m^{-1},\vert a_1 \vert_{C^{[s]+1}}, \|v\|_{H^s},\vert w\vert_{L^{\infty}}) E^s(u)^2.\nonumber
\end{align}
$\bullet$ Estimate of $\big(\Lambda^{s-2} (a_2u_{xx} ),\Lambda^0_m w^2\Lambda^su\big)$. By the Cauchy-Schwarz inequality, we have  \begin{align}
\vert\big(\Lambda^{s-2} (a_2u_{xx} ),\Lambda^0_m w^2\Lambda^su\big)\vert&\leq \frac{1}{m}\vert a_2(t,x,v,v_x)\vert_{H^{s-1}}\vert u_{xx}\vert_{H^{s-2}} \vert w\vert_{L^\infty}\vert\Lambda^s u \vert_{L^2} \nonumber\\ &\leq 
C(m^{-1},\vert a_2 \vert_{C^{[s]+1}}, \|v\|_{H^s},\vert w\vert_{L^{\infty}}) E^s(u)^2.\nonumber
\end{align}
$\bullet$Estimate of $\big(\Lambda^{s-2} (a_3u_{xxx} ),\Lambda^0_m w^2\Lambda^su\big)$.\\ We have that: $a_3u_{xxx} = \partial_x^2(a_3\partial_xu) - \partial_x^2a_3\partial_xu-2a_3\partial_x^2u$, then $\Lambda^{s-2} (a_3u_{xxx}) = \Lambda^{s-2}(\partial_x^2(a_3\partial_xu)) - \Lambda^{s-2}(\partial_x^2a_3\partial_xu)-2\Lambda^{s-2}(a_3\partial_x^2u)$.\\
Now use the identity $\Lambda^2=1-\partial_x^2$ to get that $\Lambda^{s-2}(\partial_x^2(a_3\partial_xu))=\Lambda^{s-2}\big((1-\Lambda^2)(a_3\partial_xu)\big)= \Lambda^{s-2}(a_3\partial_xu)-\Lambda^{s}(a_3\partial_xu) = \Lambda^{s-2}(a_3\partial_xu)-[\Lambda^s,a_3]\partial_xu-a_3\Lambda^s\partial_xu$, then we obtain:\\
$\big(\Lambda^{s-2} (a_3u_{xxx} ),\Lambda^0_m w^2\Lambda^su\big)=
\big(\Lambda^{s-2}(a_3\partial_xu),\Lambda^0_m w^2\Lambda^su\big)-\big([\Lambda^s,a_3]\partial_xu,\Lambda^0_m w^2\Lambda^su\big)-\big(a_3\Lambda^s\partial_xu,\Lambda^0_m w^2\Lambda^su\big)-\big(\Lambda^{s-2}(\partial_x^2a_3\partial_xu),\Lambda^0_m w^2\Lambda^su\big)-2\big(\Lambda^{s-2}(a_3\partial_x^2u),\Lambda^0_m w^2\Lambda^su\big)$.\\
By integration by parts, the third term of the last equality becomes:\\$$
\big(a_3\Lambda^s\partial_xu,\Lambda^0_m w^2\Lambda^su \big)= -\frac{1}{2}\big(\partial_x(\Lambda^0_m w^2a_3),(\Lambda^su)^2\big),$$
Now by Cauchy Shwarz we have: \begin{eqnarray*}
\vert\big(\Lambda^{s-2} (a_3u_{xxx} ),\Lambda^0_m w^2\Lambda^su\big)\vert &&\leq \frac{1}{m}\big(\| a_3\partial_xu\|_{H^{s-2}}E^s(u) +\|\partial_xa_3\|_{H^{s-1}}\|\partial_x u\|_{H^{s-1}}E^{s}(u)\\
&&+\|w^2a_3\|_{W^1,\infty}E^s(u)^2 + \|\partial_x^2a_3\partial_xu\|_{H^{s-2}}E^s(u)+ \|a_3\partial_x^2u\|_{H^{s-2}}E^s(u)\big)\\
&&\leq C(m^{-1},\|a_3\|_{C^{[s]+1}},\|v\|_{H^s},\|w\|_{W^{1,\infty}})E^s(u)^2
\\
\end{eqnarray*}
$\bullet$ Estimate of $\big([\Lambda^{s-2},a_4]\partial_x^4u,\Lambda_m^0 w^2 \Lambda^s u\big) + \big(a_4\Lambda^{s-2}\partial_x^4u,\Lambda_m^0 w^2 \Lambda^s u\big)$:\\
$a_4\Lambda^{s-2}\partial_x^4u=a_4\Lambda^{s-2}(1-\Lambda^2)\partial^2_xu= a_4(\Lambda^{s-2}-\Lambda^s)\partial^2_xu= a_4\Lambda^{s-2}\partial^2_xu - a_4\Lambda^s\partial^2_xu$, then:\\
$$\big(a_4\Lambda^{s-2}\partial_x^4u,\Lambda_m^0 w^2\Lambda^s u\big)=\big(a_4\Lambda^{s-2}\partial^2_xu,\Lambda_m^0w^2\Lambda^su\big)-\big(a_4\Lambda^s\partial^2_xu,\Lambda_m^0w^2 \Lambda^s u\big)$$
By Cauchy Shwarz, the first term of the last equality is controlled by:
$$\vert\big(a_4\Lambda^{s-2}\partial^2_xu,\Lambda_m^0w^2\Lambda^su\big)\vert \leq \frac{1}{m} \vert a_4\Lambda^{s-2}\partial^2_xu\vert_{L^2}E^s(u)\leq C(m^{-1},\vert a_4 \vert_{L^\infty})E^s(u)^2.$$
$\big(a_4\Lambda^s\partial^2_xu,\Lambda_m^0w^2 \Lambda^s u\big)=-\big(a_4\Lambda_m^0w^2,(\partial_x\Lambda^s u)^2\big)+ Q_1,$\\
where $\vert Q_1 \vert \leq C(m,s,\vert w\vert_{W^{1,\infty}},\vert \partial _x a_4\vert_{L^\infty})E^{s}(u)^2$.\\
Now, using the first order  Poisson brackets $\{\Lambda^{s-2},a_4\}_{1}= -(s-2)\partial_x(a_4)\Lambda^{s-2}\partial_x,$ see \cite{lannes'} we get:\\ $$([\Lambda^{s-2},a_4]\partial_x^4u,\Lambda_m^0w^2\Lambda^su) =(s-2)(\partial_x(a_4)\Lambda^s\partial_xu,\Lambda_m^0w^2\Lambda^su)+ Q_2,$$
Where $\vert Q_2\vert \leq C(m,s,\vert w\vert_{W^{2,\infty}},\vert a_4\vert_{H^{s+1}})E^s(u)^2$. Now, by integration by parts we have:\\$$(s-2)(\partial_x(a_4)\Lambda^s\partial_xu,\Lambda_m^0w^2\Lambda^su)= -\frac{(s-2)}{2}(\partial_x(\partial_x(a_4)\Lambda_m^0w^2)\Lambda^su,\Lambda^su)$$\\
then 
$$\vert ([\Lambda^{s-2},a_4]\partial_x^4u,\Lambda_m^0w^2\Lambda^su)\vert \leq C(m,s,\vert w\vert_{W^{2,\infty}},\vert a_4\vert_{H^{s+1}})E^s(u)^2.$$

$\bullet$ Estimate of $\big([\Lambda^{s-2},a_5]\partial_x^5u,\Lambda_m^0w^2\Lambda^su\big) + \big(a_5\Lambda^{s-2}\partial_x^5u,\Lambda_m^0w^2\Lambda^su\big):$\\
$$a_5\Lambda^{s-2}\partial_x^5u = a_5\Lambda^{s-2}(1-\Lambda^2)\partial_x^3u=a_5\Lambda^{s-2}\partial_x^3u-a_5\Lambda^s\partial_x^3u= a_5\Lambda^{s-2}\partial_xu-a_5\Lambda^{s}\partial_xu-a_5\Lambda^s\partial_x^3u$$

Then

$
\big(a_5\Lambda^{s-2}\partial_x^5u,\Lambda_m^0w^2\Lambda^su\big) = 
\big(a_5\Lambda^{s-2}\partial_xu,\Lambda_m^0w^2\Lambda^su\big)-\big(a_5\Lambda^{s}\partial_xu,\Lambda_m^0w^2\Lambda^su\big)-\big(a_5\Lambda^s\partial_x^3u,\Lambda_m^0w^2\Lambda^su\big)
$
\\
The first two terms can be easily controlled by $E^s(u)^2$ as above. Now,\\
\begin{eqnarray*}
\big(a_5\partial_x^3\Lambda^s u,\Lambda_m^0w^2\Lambda^s u\big) &&=-
\frac{1}{2}\big(\partial_x^3(a_5 \Lambda_m^0 w^2)\Lambda^s u,\Lambda^s u\big) -
\frac{3}{2}\big(\partial_x^2(w^2\Lambda_m^0 a_5)\Lambda^{s}\partial_x u,\Lambda^s u\big)\\
&&-\frac{3}{2}\big(\partial_x(\Lambda_m^0 w^2a_5)\Lambda^s u,\Lambda^{s} \partial_x^2 u\big).
\end{eqnarray*}
By integration by parts, we obtain
\begin{eqnarray*}
-\frac{3}{2}\big(\partial_x(\Lambda_m^0w^2a_5)\Lambda^s u,\Lambda^{s}\partial_x^2 u\big)=\frac{3}{2} \big(\partial_x^2(\Lambda_m^0w^2a_5)\Lambda^s u,\Lambda^{s}\partial_x u\big)
+\frac{3}{2} \big(\partial_x(a_5 \Lambda_m^0w^2),(\Lambda^{s}\partial_x u)^2\big).
\end{eqnarray*}
Now:\\
$$[\Lambda^{s-2},a_5]\partial_x^5u=\{\Lambda^{s-2},a_5\}_2\partial_x^5u+Q_3\partial_x^5u,
$$
where $\{\cdot,\cdot\}_2$ stands for the second order Poisson brackets, 
$$
\{\Lambda^{s-2},a_5\}_2=-(s-2)\partial_x(a_5)\Lambda^{s-4}\partial_x+\frac{1}{2}[(s-2)\partial_x^2(a_5)\Lambda^{s-4}-
(s-4)(s-2)\partial_x^2(a_5)\Lambda^{s-6}\partial_x^2]$$
and $Q_3$  is an operator of order $s-5$ that can be controlled by the general commutator estimates   (see \cite{lannes'}). We thus get
$$
\vert \big(Q_3\partial_x^5u,\Lambda_m^0w^2\Lambda^s u\big)\vert\leq C(m,\vert \partial_xa_5\vert_{H^{s+1}})E^s(u)^2.
$$
We now use the fact that $H^1(\R)$ is continuously embedded in $L^{\infty}(\R)$  to get 
$$
\vert \big([s\partial_x^2(a_5)\Lambda^{s-4}-
(s-4)(s-2)\partial_x^2(a_5)\Lambda^{s-6}\partial_x^2]\partial_x^5u,\Lambda_m^0w^2\Lambda^s u\big)\vert\leq C(m,s,\vert \partial_xa_5\vert_{H^{s+1}},\vert w\vert_{W^{1,\infty}})E^s(u)^2.
$$
This leads to the expression
$$
 \big([\Lambda^{s-2},a_5]\partial_x^5 u,\Lambda_m^0w^2\Lambda^su\big)=
-(s-2)\big(\partial_x (a_5) \Lambda^{s}\partial_x^2 u,\Lambda_m^0w^2\Lambda^s u\big)+Q_4,
$$
where $\vert Q_4 \vert \leq  C(m,s,\vert w\vert_{W^{1,\infty}}, \vert a_5\vert_{H^{s+1}} )E^s(u)^2$. 
Remarking now, by integration by parts
\begin{eqnarray}\label{IBP1}
-(s-2)\big(\partial_x(a_5)\Lambda^{s}\partial_x^2 u ,\Lambda_m^0w^2\Lambda^s u\big)
&&=(s-2)\big(\partial_x (\partial_x(a_5)\Lambda_m^0w^2) \Lambda^{s} \partial_xu,\Lambda^su\big)\nonumber\\
&&+(s-2)\big(\partial_x (a_3) \Lambda_m^0w^2,(\Lambda^s \partial_x u)^2\big),
\end{eqnarray}

We now choose  $w$ such that
\begin{equation}\label{eqw}
-(s-2)\big(\partial_x (a_5)\Lambda_m^0 w^2,(\Lambda^s \partial_x u)^2\big)+\frac{3}{2} \big(\partial_x(a_5\Lambda_m^0w^2),(\Lambda^{s}\partial_x u)^2\big) + \big(a_4\Lambda_m^0w^2,(\partial_x\Lambda^s u)^2\big)=0;
\end{equation}
therefore, if we take $w=(\Lambda_m^0)^{-1}\Big(\vert a_5\vert^ {\big({\frac{2s-7}{6}\big)}}\displaystyle{\exp(-\frac{1}{3}\int_0^x\frac{a_4}{a_5}dy)}\Big)$ we  easily  obtain (\ref{eqw}). Finally, one has
\begin{eqnarray*}
 &&\big([\Lambda^{s-2},a_5]\partial_x^5 u,\Lambda_m^0w^2\Lambda^s  u\big)+\big(a_5\partial_x^5\Lambda^{s-2} u,\Lambda_m^0w^2\Lambda^s u\big)\\
&&\hspace*{3cm} =Q_4+(s-2)\big(\partial_x (\partial_x(a_5)\Lambda_m^0w^2) \Lambda^s\partial_x u,\Lambda^s u\big)
-\frac{1}{2}\big(\partial_x^3(a_5\Lambda_m^0w^2)\Lambda^s u,\Lambda^s u\big)\\&&\hspace*{3cm} \qquad-
\frac{3}{2}\big(\partial_x^2(a_5\Lambda_m^0w^2)\Lambda^{s}\partial_x u,\Lambda^s u\big)+
\frac{3}{2} \big(\partial_x^2(a_5\Lambda^0_mw^2)\Lambda^{s}\partial_x u,\Lambda^su\big);
\end{eqnarray*}
therefore,
$$
\vert\big([\Lambda^{s-2},a_5]\partial_x^5 u,\Lambda_m^0w^2\Lambda^s  u\big)+\big(a_5\partial_x^5\Lambda^{s-2} u,\Lambda_m^0w^2\Lambda^s u\big)\vert \leq C(s,m,\vert \partial_xa_5\vert_{H^{s+1}})
E^s(u)^2.
$$
$\bullet$ Estimate of $\big(w_t\Lambda^{s-2}u ,\Lambda_m^0w\Lambda^su\big)$: Using the Cauchy-Schwarz inequality we obtain
$$
\vert \big(w_t\Lambda^su ,w\Lambda^su\big)\vert \leq  C(m,\vert w_t\vert_{L^{\infty}},\vert w\vert_{L^{\infty}})
E^s(u)^2.
$$
Gathering the information provided by the above estimates, since one has $$\vert \big(\Lambda^{s-2}f ,\Lambda_m^0w^2\Lambda^su\big)\vert \leq 
\frac{1}{m}E^s(f)E^s(u).
$$
If we assemble the previous estimates  and using Gronwall's lemma we obtain the following estimate:\\
$$
e^{\lambda t}\partial_t (e^{-\lambda t}E^s(u)^2) 
	\leq  \big(C(E^s(v))-\lambda\big)E^s(u)^2+2 E^s(f)E^s(u).
$$
Taking $\lambda=\lambda_T$ large enough (how large depends on 
$\sup_{t\in [0,T]}C(E^s(v(t))$
for the first term of the right hand side of the above inequality to be negative for all $t\in [0,T]$, we deduce that
$$
E^s(u(t)) 
	\leq e^{\lambda_T t}E^s(u^0)
	+2 \int_0^t e^{\lambda_T (t-t')}
	E^s(f(t'))dt'.
$$
\subsection{Proof of the theorem:}

Thanks to this energy estimate, we classically conclude 
(see e.g. \cite{AG}) the existence of a time 
$$
	T^*=T^*(E^s(u^0))>0,
$$ 
and a unique solution $u\in C([0,T^*];H^{s}(\R))\cap
	C^1([0,T^*];H^{s-3}(\R))$ to
(\ref{hch}) as the limit of the iterative scheme
$$
	u_0=u^0,\quad\mbox{ and }\quad
	\forall n\in\N, \quad
	\left\lbrace\begin{array}{l}
	{\mathcal L}(u^n,\partial)u^{n+1}=f,\\
	u^{n+1}_{\vert_{t=0}}=u^0.
		    \end{array}\right.
$$

\providecommand{\href}[2]{#2}

\end{document}